\documentstyle[11pt,leqno]{article}
\newtheorem{theorem}{{\sc Theorem}}
\newcommand{\bt}{\begin{theorem}}
\newcommand{\et}{\end{theorem}}
\newcommand{\newsection}[1]{\setcounter{equation}{0} \setcounter{theorem}{0}
\section{#1}}

\newcommand{\NI}{\noindent}
\newcommand{\bea}{\begin{eqnarray}}
\newcommand{\eea}{\end{eqnarray}}

\def \spec#1 {\mathop{#1}}

\def \b #1 {\bf #1}
\newcommand {\nnb }{\nonumber}

\newcommand {\CC}{\centerline}

\newcommand{\clf}{{\cal F}}

\newcommand{\ity}{\infty}
\newcommand{\raro}{\rightarrow}

\newcommand{\vsp}{\vskip 1em}

\newcommand{\ve}{\varepsilon}

\newcommand{\be}{\begin{equation}}
\newcommand{\ee}{\end{equation}}
\newcommand{\ben}{\begin{eqnarray*}}
\newcommand{\een}{\end{eqnarray*}}

\begin{document}

\sloppy
\CC {\Large {\bf  Minimum $L_1$-norm Estimation for Fractional}}
\CC{\Large {\bf Ornstein-Uhlenbeck Process Driven by a Gaussian Process}} 
\vsp
\CC {\bf B.L.S. Prakasa Rao   }
\CC{\bf CR Rao Advanced Institute of Mathematics, Statistics }
\CC{\bf and Computer Science, Hyderabad, India}
\vsp
\NI{\bf Abstract:} We investigate the asymptotic properties of the
minimum $L_1$-norm estimator of the drift parameter for fractional
Ornstein-Uhlenbeck type process driven by a general Gausssian process.
\vsp
\newsection{Introduction}

Diffusion processes and diffusion type processes satisfying stochastic
differential equations driven by Wiener processes are used for stochastic
modeling in a wide variety of sciences such as population genetics,
economic processes, signal processing as well as for modeling sunspot
activity and more recently in mathematical finance. Statistical inference
for diffusion type processes satisfying stochastic  differential equations
driven by Wiener processes have been studied earlier and a comprehensive
survey of various methods is given in Prakasa Rao (1999). There has been
a recent interest to study similar problems for stochastic processes
driven by a fractional Brownian motion to model processes involving long
range dependence (cf. Prakasa Rao (2010)). Le Breton (1998) studied parameter estimation and
filtering in a simple linear model  driven by a fractional Brownian
motion. Kleptsyna and Le Breton (2002) studied
parameter estimation problems for fractional Ornstein-Uhlenbeck process.
The fractional Ornstein-Uhlenbeck process is a fractional analogue of the Ornstein-Uhlenbeck process, that is,
a continuous time first order autoregressive process $X=\{X_t, t \geq 0\}$
which is the solution of a one-dimensional homogeneous linear stochastic
differential equation driven by a fractional Brownian motion (fBm) $W^H=
\{W_t^H, t \geq 0\}$ with Hurst parameter $H.$ Such a process
is the unique Gaussian  process satisfying the linear integral equation

\be
X_t= x_0+ \theta \int_0^tX_s ds + \sigma W_t^H, t \geq 0.
\ee
They investigate the problem of estimation of the parameters $\theta$ and
$\sigma^2$ based on the observation $\{X_s, 0 \leq s \leq T\}$ and study the asymptotic behaviour of these estimators
as $T \raro \ity.$

In spite of the fact that maximum likelihood estimators (MLE) are consistent
and asymptotically normal and also asymptotically efficient in general,
they have some short comings at the same time. Their calculation is often
cumbersome as the expression for MLE involve stochastic integrals at times which
need good approximations for computational purposes. Further more MLE are
not robust in the sense that a slight perturbation in the noise component
will change the properties of MLE substantially. In order to circumvent
such problems, the minimum distance approach is proposed. Properties of
the minimum distance estimators (MDE) were discussed in Millar (1984) in a
general frame work. Kutoyants and Pilibossian (2000) studied the problem of minimum $L_1$-norm estimation for the Ornstein-Uhlenbeck process. Prakasa Rao (2004) investigated the problem of minimum $L_1$-norm estimation for 
the fractional Ornstein-Uhlenbeck process driven by a fractional Brownian motion.

Our aim in this paper is to obtain the minimum $L_1$-norm estimator of the
drift parameter of a Ornstein-Uhlenbeck  process driven by general Gaussian processes and
investigate the asymptotic properties of such estimators. El Machkouri et al. (2015), Chen and Zhou (2020) and Lu (2022) study parameter estimation for an Ornstein-Uhlenbeck process  driven by a general Gaussian process.

\newsection{Minimum $L_1$-norm Estimation}
Let $(\Omega, \clf, (\clf_t), P) $ be a stochastic basis satisfying the
usual conditions and the processes discussed in the following are
$(\clf_t)$-adapted. Further the natural filtration of a process is
understood as the $P$-completion of the filtration generated by this
process. We consider centered a Gaussian process $G\equiv\{G_t,0\leq t \leq 1\}.$
\vsp
Let us consider a stochastic process $\{X_t, t \in [0,1]\}$ defined by the 
stochastic integral equation 
\be
X_t=x_0 +  \theta \int_0^t X(s) ds  + \ve G_t, 0 \leq t \leq 1,
\ee
where $\theta $ is an unknown  drift parameters respectively. For
convenience, we write the above integral equation in the form of a
stochastic differential equation   
\be
dX_t =\theta X(t) dt + \ve  dG_t, X_0=x_0, 0 \leq t \leq 1,
\ee
driven by the Gaussian process $G.$ The class of Gaussian processes $G$ includes fractional Brownian motion, sub-fractional Brownian motion and bifractional Brownian motion.
\vsp
We now consider the problem of estimation of the parameter $\theta $ based
on the observation of fractional Ornstein-Uhlenbeck process $ X=
\{X_t , 0 \leq t \leq 1\}$ satisfying the stochastic differential equation 
\be
dX_t=  \theta  X(t) dt  + \ve dG_t, X_0=x_0, 0 \leq t \leq 1
\ee
where $\theta \in \Theta \subset R$ and study its
asymptotic properties as $\ve \raro 0.$ 
\vsp
Let $x_t(\theta)$ be the solution of the above differential equation with
$\ve=0 .$ It is obvious that 
\be
x_t(\theta)= x_0e^{\theta t}, 0 \leq t \leq 1.
\ee
Let 
\be
S_1(\theta)= \int_0^1|X_t-x_t(\theta)|dt
\ee
\vsp
We define $\theta_\ve^*$ to be a {\it minimum $L_1$-norm estimator} if 
there exists a measurable selection $\theta_\ve^*$ such that  
\be
S_1(\theta_\ve^*)= \inf_{\theta \in \Theta}S_1(\theta). 
\ee 
Conditions for the existence of a measurable selection are given in Lemma
3.1.2 in Prakasa Rao (1987). We assume that there exists a measurable
selection $\theta_\ve^*$ satisfying the above equation. An alternate way
of defining the estimator $\theta_\ve^*$ is by the relation   
\be
\theta_\ve^*= \arg \inf_{\theta \in \Theta} \int_0^1|X_t-x_t(\theta)|dt.
\ee
\vsp
We will now present some maximal inequalities for Gaussian processes out of which one will be used in the sequel.
\vsp
\NI{\bf Some Maximal Inequalities for Gaussian processes:}
\vsp
Let $G_1^*= \sup_{0\leq t \leq 1}|G_t|.$ If $G$ is a fractional Brownian motion or a sub-fractional Brownian motion, maximal inequalities are known and are reviewed in Prakasa Rao (2014, 2017, 2020) following self-similarity for such processes. Maximal inequalities for general Gaussian processes are surveyed in Li and Shao (2001). We now present some maximal inequalities for general Gaussian processes $G.$ Let
$$d_G(s,t)= [E(G_t-G_s)^2]^{1/2}, 0\leq s,t\leq 1,$$
$$\rho(\ve)= \sup_{t,s \in [0,1], |s-t|\leq \ve}d_G(s,t)$$
and
$$Q(\delta)=\int_0^\ity\rho(\delta e^{-y^2})dy.$$
Suppose that the function $\rho(.)$ is strictly increasing. Let $\sigma^2= \sup_{0\leq t \leq 1}E(G_t^2)$ and $Q^{-1}(.)$ denote the inverse of the function $Q(.)$. The following result is due to Berman (1985).
\vsp
\NI {\bf Theorem 2.1:} Under the conditions stated above, given a Gaussian process $G$ defined on the interval $[0,1],$ there exists a constant $C$ such that
$$P(\sup_{0\leq t \leq 1}|G_t|> \ve) \leq C(Q^{-1}(1/\ve))^{-1}\exp(-\frac{\ve^2}{2\sigma^2}).$$
for all $\epsilon >0.$
\vsp
Consider a Gaussian process $G$ on the interval $[0,1]$ and let $D= \sup_{0\leq s,t\leq 1}d_G(s,t).$ Let $N(\ve,d_G,1)$ denote the minimum  number of closed intervals  of length $2\ve >0,$ needed to cover the interval $[0,1]$ with mid-points in $[0,1]$ in the semi-metric $d_G.$ We say that the Gaussian process $G$ is pairwise non-degenerate if for every $\epsilon >0$ there exist closed intervals  of length $2\ve >0,$ in the semi-metric $d_G$ needed to cover the Interval $[0,1].$ If the Gaussian process is pairwise non-degenerate, then the semi-metric $d_G$ is a metric. The following result is due to Marcus and Rosen (2006).
\vsp
\NI{\bf Theorem 2.2:} Suppose $G$ is a Gaussian process which is pairwise non-degenerate and continuous in mean square. If 
$$\int_0^D\sqrt{\log N(\epsilon,d_G,1)}d\epsilon <\ity,$$
then the process $G$ has continuous sample paths and there exists a positive constant $C$ such that 
$$E[\sup_{0\leq s \leq 1}|G(s)|]\leq C \int_0^{D/2}\sqrt{\log N(\epsilon,d_G,1)}d\epsilon <\ity.$$
\vsp
The following result is due to Borovkov et al. (2017).
 \vsp
\NI{\bf Theorem 2.3:} Suppose $G$ is a Gaussian process on the interval $[0,1]$ which is pairwise non-degenerate and there exists constants $0\leq C_1\leq C_2$ and $H\in (0,1)$ such that 
$$C_1|t-s|^H \leq d_G(s,t)\leq C_2|t-s|^H, t, s \in [0,1].$$
Then
$$\frac{C_1}{5 \sqrt{H}} \leq E[\sup_{0\leq t \leq 1}G_t] \leq \frac{(16.3) C_2}{\sqrt{H}}.$$
\vsp

Another useful maximal inequality for Gaussian process is the following result due to Nourdin (2012). 

\NI{\bf Theorem 2.4:} Suppose $G$ is a centered and continuous Gaussian process on the interval $[0,1]$.  Let $\sigma^2= \sup_{0\leq t \leq 1}E[G_t^2].$ Suppose that $0<\sigma^2 <\ity.$ Then $m= E[\sup_{0\leq t \leq 1}G_t]$ is finite and for all $x >m,$
$$P(\sup_{0\leq t \leq 1}G_t\geq x)\leq \exp(-\frac{(x-m)^2}{2\sigma^2}).$$
\vsp
From the fact that $G$ is centered Gaussian process, it follows that
$$P(\sup_{0\leq t \leq 1}|G_t|\geq x)\leq 2 \; \exp(-\frac{(x-m)^2}{2\sigma^2})$$
under the conditions stated in Theorem 2.4.
\vsp
\newsection{\bf Consistency of the estimator:}

Let $\theta_0$ denote the true parameter. For any $\delta >0,$ define
\be
g(\delta)= \inf_{|\theta-\theta_0|>\delta}\int_0^1|x_t(\theta)-x_t(\theta_0)|dt.  
\ee
Note that $g(\delta ) >0$ for any $\delta >0.$ 
\vsp
\NI{\bf Theorem 3.1:} Suppose $G$ is a centered and continuous Gaussian process on the interval $[0,1]$. Let $\sigma^2= \sup_{0\leq t \leq 1}E[G_t^2] $ and $m=E[\sup_{0\leq t \leq 1}G_t].$ Then there exists a positive constant $C$ such that for every  $\delta >0,$
$$P_{\theta_0}^{(\ve)}\{|\theta_\ve^*-\theta_0|> \delta \} = O(e^{-C[g(\delta)]^2\ve^{-2}}).$$
\vsp
\NI {\bf Proof:} Let $||.||$ denote the $L_1$-norm. Then
\bea
P_{\theta_0}^{(\ve)}\{|\theta_\ve^*-\theta_0|> \delta\} & = &
P_{\theta_0}^{(\ve)}\{\inf_{|\theta -\theta_0|\leq  \delta}
||X - x(\theta)|| > \inf_{|\theta -\theta_0| >  \delta} ||X-x(\theta)|| \} \\ \nnb
& \leq &  P_{\theta_0}^{(\ve)}\{\inf_{|\theta -\theta_0|\leq  \delta}
(||X - x(\theta_0)|| + || x(\theta)-x(\theta_0)||) \\ \nnb
& & \;\;\;\;\; > \inf_{|\theta -\theta_0| >  \delta}
(||x(\theta) - x(\theta_0)|| - || X-x(\theta_0)||)\} \\ \nnb
& = & P_{\theta_0}^{(\ve)}\{
2||X - x(\theta_0)|| > \inf_{|\theta -\theta_0| >  \delta}
||x(\theta)-x(\theta_0)|| \} \\ \nnb
& = & P_{\theta_0}^{(\ve)}\{
||X - x(\theta_0)|| > \frac{1}{2} g(\delta)\}.
\eea
\vsp
Since the process $X$ satisfies the stochastic differential equation
(2.1), it follows that 
\bea
X_t -x_t(\theta_0) &=& x_0+ \theta_0 \int_0^t X_s ds + \ve G_t -x_t(\theta_0) \\ \nnb
& = & \theta_0\int_0^t(X_s-x_s(\theta_0))ds + \ve G_t
\eea
since $x_t(\theta)= x_0 e^{\theta t}.$ Let $U_t= X_t-x_t(\theta_0).$ Then
it follows from the equation given above that
\be
U_t= \theta_0 \int_0^t U_s\;\;ds + \ve G_t.
\ee
Let $V_t=|U_t|=|X_t-x_t(\theta_0)|.$ The relation given above implies that
\be
V_t= |X_t-x_t(\theta_0)| \leq |\theta_0| \int_0^t V_s ds + \ve |G_t|.
\ee
Applying the Gronwall-Bellman Lemma, it follows that
\be
\sup_{0\leq t \leq 1} |V_t| \leq \ve e^{|\theta_0 |} \sup_{0 \leq t \leq
1}|G_t|.
\ee
Hence
\bea
P_{\theta_0}^{(\ve)}[||X - x(\theta_0)|| > \frac{1}{2} g(\delta)]
&\leq & P[\sup_{0 \leq t \leq 1}|G_t| >  \frac{e^{-|\theta_0|}g(\delta)}{2 \ve }] \\ \nnb
& = & P[G_1^* >  \frac{e^{-|\theta_0|}g(\delta)}{2 \ve }].
\eea
Let $m= E[\sup_{0\leq t \leq 1} G(t)].$ Applying the  maximal inequalities for  Gaussian processes given in Theorem 2.4, we get that, for fixed  $\delta >0,$ we can choose $\epsilon$ sufficiently small so that $\frac{e^{-|\theta_0|} g(\delta)}{2\epsilon}>m.$ For such $\ve,$
\bea
P_{\theta_0}^{(\ve)}[|\theta_\ve^*-\theta_0|> \delta] & \leq & 2 \exp(-\frac{((e^{-|\theta_0|}g(\delta)/2\epsilon)-m)^2}{2\sigma^2})\\ \nnb
& = & O(e^{-C[g(\delta)]^2\ve^{-2}})
\eea
for some positive constant $C$ independent of $\ve.$
\vsp
\NI{\bf Remarks: } As a consequence of the result obtained above, it follows that 
$$P_{\theta_0}^{(\ve)}\{|\theta_\ve^*-\theta_0|> \delta\} \raro 0\;\; \mbox{as}\;\; \epsilon \raro 0$$
for every $\delta >0.$ Hence the minimum norm $L_1$-estimator $\theta_\ve^*$ is weakly consistent for estimating the parameter $\theta_0.$

\newsection{\bf Asymptotic distribution of the estimator:}

We will now study the asymptotic distribution if any of the estimator
$\theta_\ve^*$ after suitable scaling. It can be checked that
\be
X_t= e^{\theta_0t}\{x_0+\int_0^te^{-\theta_0s} \ve dG_s\}
\ee
or equivalently
\be
X_t-x_t(\theta_0)= \ve e^{\theta_0t} \int_0^te^{-\theta_0 s}dG_s.
\ee
Let 
\be
Y_t= e^{\theta_0t} \int_0^te^{-\theta_0 s}dG_s.
\ee
Note that $\{Y_t, 0 \leq t \leq 1\}$ is a Gaussian process and can be
interpreted as the "derivative" of the process $\{X_t, 0 \leq t \leq 1\}$
with respect to $\ve.$ We obtain that, $P$-a.s., 
\be
Y_t e^{-\theta_0 t}= \int_0^te^{-\theta_0 s}dG_s 
\ee
The integral with respect to the process $G$ is interpreted as Young integral (cf. El Machkouri et al. (2015)). In particular it
follows that the random variable $Y_t e^{-\theta_0t}$ and hence $Y_t$ has
the normal distribution with mean zero and further more, for any $0 \leq t,s\leq 1,$ 
\bea
Cov(Y_t,Y_s) &=& e^{\theta_0 t + \theta_0 s} E[\int_0^te^{-\theta_0 u}dG_u   \int_0^{s}e^{-\theta_0 v}dG_v] \\ \nnb
& = & e^{\theta_0 t + \theta_0 s} \int_0^t \int_0^s e^{-\theta_0(u+v)} dudv \\ \nnb
& = &  R(t,s)\;\;\mbox{(say)}. 
\eea
In particular
\be
Var(Y_t)= R(t,t).
\ee
Observe that $\{Y_t, 0 \leq t \leq 1\}$ is a zero mean Gaussian process with $Cov(Y_t,Y_s)= R(t,s).$ Let
\be
\zeta= \arg \inf_{-\ity < u < \ity} \int_0^1|Y_t- ut x_0e^{\theta_0 t}|dt.
\ee
\vsp
\NI{\bf Theorem 4.1:} As $\ve \raro 0,$ the random variable
$\ve^{-1}(\theta_\ve^*-\theta_0)$ converges in probability to a
random variable whose probability distribution is the same as that of
$\zeta$ under $P_{\theta_0}.$
\vsp
\NI{\bf Proof:} Let $x^\prime_t(\theta)=x_0 t e^{\theta t}$ and let
\be
Z_{\ve}(u)= ||Y-\ve^{-1}(x(\theta_0+\ve u)-x(\theta_0))||
\ee
and
\be
Z_0(u)= ||Y-u x^\prime(\theta_0)||.
\ee
Further more, let
\be
A_{\ve}= \{\omega: |\theta_\ve^*-\theta_0| < \delta_{\ve}\}, \delta_\ve=
\ve^\tau, \tau \in (\frac{1}{2},1), L_\ve= \ve^{\tau-1}.
\ee
Observe that the random variable $u_\ve^*=
\ve^{-1}(\theta_\ve^*-\theta_0)$ satisfies the equation 
\be
Z_{\ve}(u_\ve^*)= \inf_{|u|<L_\ve} Z_\ve(u), \omega \in A_\ve.
\ee
Define
\be
\zeta_\ve= \arg \inf_{|u|<L_\ve}Z_0(u).
\ee
Observe that, with probability one,
\bea
\sup_{|u|< L_\ve}|Z_\ve(u)-Z_0(u)| & = & | || Y-u x^\prime
(\theta_0)-\frac{1}{2}\ve u^2 x^{\prime \prime}(\tilde \theta)|| - ||Y-u
x^\prime( \theta_0)||| \\ \nnb
& \leq & \frac{\ve}{2}L_\ve^2\sup_{|\theta-\theta_0<\delta_\ve}\int_0^T|x^{\prime\prime}
(\theta)|dt \\ \nnb
& \leq & C \ve^{2 \tau -1}.
\eea
Here $\tilde \theta = \theta_0+ \alpha(\theta-\theta_0)$ for some $\alpha
\in (0,1).$ Note that the last term in the above inequality tends to zero
as $\ve \raro 0.$ Further more the process $\{Z_0(u), -\ity < u <\ity
\}$ has a unique minimum $u^*$ with probability one. This follows from the
arguments given in Theorem 2 of Kutoyants and Pilibossian (1994). In addition, we can choose the interval $[-L,L]$ such that
\be
P_{\theta_0}^{(\ve)} \{u_\ve^* \in (-L,L)\} \geq 1- \beta (g(L))^{-1}  
\ee
and
\be
P\{u^* \in (-L,L)\} \geq 1- \beta (g(L))^{-1}
\ee
where $\beta >0.$ Note that $g(L)$ increases as $L$ increases. The
processes $Z_\ve(u), u \in [-L,L]$ and $Z_0(u), u \in [-L,L]$ satisfy the
Lipschitz conditions and $Z_\ve(u)$ converges uniformly to $Z_0(u)$ over $u
\in [-L,L].$ Hence the minimizer of $Z_\ve(.)$ converges to the minimizer
of $Z_0(u).$ This completes the proof.  
\vsp
\NI{\bf Remarks :} We have seen earlier that the process $\{Y_t, 0 \leq t \leq T\}$ is a zero mean Gaussian process with the covariance function $Cov(Y_t,Y_s)= R(t,s)$ for $0,t,s \leq 1.$ Recall that
\be 
\zeta= \arg \inf_{-\ity < u < \ity} \int_0^1|Y_t- utx_0 e^{\theta_0 t}|dt.
\ee
It is not clear what the distribution of the random variable $\zeta $ is. It depends on the Gaussian process $G.$ Observe that for
every $u,$ the integrand in the above integral is the absolute value of a Gaussian process $\{J_t, 0 \leq t \leq 1 \}$  with the mean function
$E(J_t)= -utx_0e^{\theta_0t}$ and the covariance function $Cov(J_t,J_s)= R(t,s)$ for $0 \leq s,t \leq 1.$ It is easy to extend the results to any Gaussian process defined on any interval $[0,T]$ for any fixed $T>0.$
\vsp
\NI{\bf Acknowledgment:} This work was supported by the INSA Senior Scientist fellowship at the CR Rao Advanced Institute of Mathematics, Statistics and Computer Science, Hyderabad, India.
\vsp
\NI {\bf References}
\begin{description}

\item Berman, S.M. (1985) An asymptotic bound for the tail of the distribution of the maximum of a Gaussian process, {\it Ann. Inst. H. Poincare Probab. Statist.}, {\bf 21}, 47-57.

\item Borovkov, K., Mishura, Y., and Novikov, A. (2017) Bounds for expected maximum of Gaussian processes and their discrete approximations, {\it Stochastics}, {\bf 89}, 21-37.

\item Chen, Y., and Zhou, H. (2020) Parameter  estimation for an Ornstein-Uhlenbeck process driven by a general Gaussian noise, arxiv:2002.09641v1 [math.PR] 22 Feb 2020.

\item El Machkouri, M., Es-Sebaiy, K., and Ouknine, Y. (2015) Parameter estimation for the non-ergodic Ornstein-Uhlenbeck processes driven by Gaussian process, arXiv:1507.00802v1 [math.PR] 3 July 2015.

\item Kleptsyna, M.L. and Le Breton, A. (2002) Statistical analysis of the fractional Ornstein-Uhlenbeck type process, {\it Statist. Infer.
Stoch. Proces.} {\bf 5}, 229-248. 

\item  Kutoyants, Yu. and Pilibossian, P. (1994) On minimum $L_1$-norm estimate of the parameter of the Ornstein-Uhlenbeck process, {\it Statist. Probab. Lett.}, {\bf 20}, 117-123.

\item  Le Breton, A. (1998) Filtering and parameter estimation in a simple linear model driven by a fractional Brownian motion, {\it Statist. Probab. Lett.} {\bf 38}, 263-274.

\item Lu, Y. (2022) Parameter estimation of non-ergodic Ornstein-Uhlenbeck processes driven by general Gaussian processes, arXiv:2207:13355v1 [math.ST] 27 Jul 2022. 

\item Marcus, M.B. and Rosen, J. (2006) {\it Markov processes, Gaussian processes and Local Times}, Cambridge Studies in Advanced Mathematics, Vol. 100, Cambridge University Press, Cambridge. 

\item  Millar, P.W. (1984) A general approach to the optimality of the minimum distance estimators, {\it Trans. Amer. Math. Soc.} {\bf 286}, 
249-272.

\item Nourdin, I. (2012) {\it Selected Aspects of Fractional Brownian Motion}, Bocconi and Springer Series, Bocconi University Press, Milan.

\item  Prakasa Rao, B.L.S. (1987) {\it Asymptotic Theory of Statistical Inference}, Wiley, New York.

\item  Prakasa Rao, B.L.S. (1999)  {\it Statistical Inference for Diffusion Type Processes}, Arnold, London and Oxford University Press, New
York. 

\item  Prakasa Rao, B.L.S. (2004) Minimum $L_1$-norm estimation for fractional Ornstein-Uhlenbeck type process, {\it Theor. Probability and Math. Statist.}, {\bf 71} (2004) 160-168. 

\item Prakasa Rao, B.L.S. (2010) {\it Statistical Inference for Fractional Diffusion Processes}, Wiley, London.

\item Prakasa Rao, B.L.S. (2014) Maximal inequalities for fractional Brownian motion: An overview, {\it Stochastic Analysis and Applications}, {\bf 32}, 450-479.

\item Prakasa Rao, B.L.S. (2017) On some maximal and integral inequalities for sub-fractional Brownian motion, {\it Stochastic Analysis and Applications}, {\bf 35},  279-287.

\item Prakasa Rao, B.L.S. (2020) More on maximal inequalities for sub-fractional Brownian motion, {\it Stochastic Analysis and Applications}, {\bf 38}, 238-247.
\end{description}
\vsp
\NI {CR Rao Advanced Institute of Mathematics, Statistics and Computer Science, Hyderabad, India.}\\ 
\NI{e-mail: blsprao@gmail.com}
\end{document}